\documentclass[12pt]{article}
\usepackage{amsthm}
\usepackage{amssymb}
\usepackage{amsmath}

\newtheorem{thm}{Theorem}

\newtheorem{defi}{Definition}
\newtheorem{cor}{Corollary}
\newtheorem{rem}{Remark}

\newtheorem{lem}{Lemma}

\newcommand{\idem}{\operatorname{idem}}

\begin{document}
\title{A connection formula of a divergent bilateral basic hypergeometric function}
\author{Takeshi MORITA\thanks{Graduate School of Information Science and Technology, Osaka University, 
1-1  Machikaneyama-machi, Toyonaka, 560-0043, Japan.} }
\date{}
\maketitle
\begin{abstract}
We give the new connection formula for the divergent bilateral basic hypergeometric series ${}_2\psi_2(a_1,a_2;b_1;q,x)$ by the using of the $q$-Borel-Laplace resummation method and Slater's formula. The connection coefficients are given by elliptic functions.  
\end{abstract}

\section{Introduction}
In this paper, we show a connection formula for a divergent \textit{bilateral} basic hypergeometric function 
\begin{equation}
{}_2\psi_1(a_1,a_2;b_1;q,x):=\sum_{n\in\mathbb{Z}}\frac{(a_1;q)_n(a_2;q)_n}{(b_1;q)_n(q;q)_n}\left\{(-1)^nq^{\frac{n(n-1)}{2}}\right\}^{-1}x^n.
\label{div}\end{equation}
Here, $(a;q)_n, (n\in\mathbb{Z})$ is the $q$-shifted factorial (see \cite{GR} for more details). We assume that $q\in\mathbb{C}^*$ satisfies $0<|q|<1$. The function \eqref{div} satisfies the second order linear $q$-difference equation
\begin{equation}
\left(\frac{b_1}{q^2}-a_1a_2x\right)u(q^2x)-\left\{\frac{1}{q}-(a_1+a_2)x\right\}u(qx)-xu(x)=0.\label{diveq}
\end{equation}
The equation \eqref{diveq} also has the \textit{unilateral} solutions around infinity:

\begin{align*}v_1(x)&=\frac{\theta (a_1x)}{\theta (x)}\sum_{n\ge 0}\frac{(qa_1/b_1;q)_n (b_1/a_1a_2x)^n}{(qa_1/a_2;q)_n(q;q)_n},\\
v_2(x)&=\frac{\theta (a_2x)}{\theta (x)}\sum_{n\ge 0}\frac{(qa_2/b_1;q)_n(b_1/a_1a_2x)^n}{(qa_2/a_1;q)_n(q;q)_n},
\end{align*}
provided that the function $\theta (x):=\sum_{n\in\mathbb{Z}}q^{n(n-1)/2}x^n$, $\forall x\in\mathbb{C}^*$ is the theta function of Jacobi.
\noindent
The aim of this paper is to give the connection formula between $v_1(x)$, $v_2(x)$ and the divergent series \eqref{div} as follows:

\bigskip
\noindent
\textbf{Theorem.}
For any $x\in\mathbb{C}^*\setminus -\lambda q^{\mathbb{Z}}$, we have 

\begin{align*}
&\left(\mathcal{L}_{q,\lambda}^+\circ\mathcal{B}^+_q{}_2\psi_1(a_1,a_2;b_1;q,x)\right)(x)\\
=&\frac{(1/a_2,qa_1/a_2,b_1/a_1,q;q)_\infty}{(b_1,q/a_1,a_1/a_2,qa_2/a_1;q)_\infty}\frac{\theta (a_1\lambda /q)}{\theta (\lambda /q)}
\frac{\theta (a_1qx/\lambda )}{\theta (qx/\lambda )}\frac{\theta (x)}{\theta (a_1x)}v_1(x)\notag\\
+&\frac{(1/a_1,qa_2/a_1,b_1/a_2,q;q)_\infty}{(b_1,q/a_2,a_2/a_1,qa_1/a_2;q)_\infty}\frac{\theta (a_2\lambda /q)}{\theta (\lambda /q)}
\frac{\theta (a_2qx/\lambda )}{\theta (qx/\lambda )}\frac{\theta (x)}{\theta (a_2x)}v_2(x).\notag
\end{align*}
Here, $\mathcal{B}_q^+$ and $\mathcal{L}_{q,\lambda}^+$ are \textit{the $q$-Borel-Laplace transformations}(see section two). We remark that the $q$-elliptic functions (with the new parameter $\lambda$) appear in the connection coefficients. 

\bigskip

At first, we review the connection problems on the linear $q$-difference equations. Connection problems on the linear $q$-difference equations with regular singular points were studied by G.~D.~Birkhoff \cite{Birkhoff} in 1914. Connection formulae for the second order linear $q$-difference equations are given by the following matrix form:

\[\begin{pmatrix}
u_1(x)\\
u_2(x)
\end{pmatrix}
=
\begin{pmatrix}
C_{11}(x)&C_{12}(x)\\
C_{21}(x)&C_{22}(x)
\end{pmatrix}
\begin{pmatrix}
v_1(x)\\
v_2(x)
\end{pmatrix}.\]
The pair ($u_1(x), u_2(x)$) is a fundamental system of (unilateral) solutions around the origin and the pair $(v_1(x), v_2(x))$ is a fundamental system of solutions around the infinity. The connection coefficients $C_{jk}$ $(1\le j,k\le 2)$ are elliptic functions:
\[\sigma_qC_{jk}(x)=C_{jk}(x),\quad C_{jk}(e^{2\pi i}x)=C_{jk}(x),\]
namely, $q$-periodic and unique valued functions.

The first example of the connection formula was given by G.~N.~Watson \cite{W} in 1910. Watson gave the connection formula for Heine's basic hypergeometric series 
\[{}_2\varphi_1(a,b;c;q,x):=\sum_{n\ge 0}\frac{(a,b;q)_n}{(c;q)_n(q;q)_n}x^n\]
around the origin and around the infinity \cite[page 117]{GR}. 
Heine's ${}_2\varphi_1(a,b;c;q,x)$ satisfies the $q$-difference equation 
\begin{equation}
\left[(c-abqx)\sigma_q^2-\left\{(c+q)-(a+b)qx\right\}\sigma_q+q(1-x)\right]u(x)=0.\label{heineeq}
\end{equation}
The equation \eqref{heineeq} also has a fundamental system of solutions around the infinity:
\[y_{\infty}^{(a,b)}(x)=\frac{\theta (ax)}{\theta (x)}{}_2\varphi_1\left(a,\frac{aq}{c};\frac{aq}{b};q,\frac{cq}{abx}\right)\]
and
\[y_{\infty}^{(b,a)}(x)=\frac{\theta (bx)}{\theta (x)}{}_2\varphi_1\left(b,\frac{bq}{c};\frac{bq}{a};q,\frac{cq}{abx}\right).\]
Watson's connection formula for ${}_2\varphi_1(a,b;c;q,x)$ is given by
\begin{align}\label{wato}
{}_2 \varphi_1\left(a,b;c;q;x \right)&= 
\frac{(b,c/a;q)_\infty \theta (-ax)_\infty }{(c, b/a;q)_\infty \theta (-x)_\infty }\frac{\theta (x)}{\theta (ax)} 
y_{\infty}^{(a,b)}(x) \nonumber \\
&+\frac{(a,c/b;q)_\infty \theta (-bx)_\infty }{(c, a/b;q)_\infty \theta (-x)_\infty } \frac{\theta (x)}{\theta (bx)} 
y_{\infty}^{(b,a)}(x).\notag
\end{align}
Here, connection coefficients are $q$-elliptic functions.

But connection formulae for $q$-difference equations with irregular singular points had not known for a long time. 
The irregularity of $q$-difference equations are studied using the Newton polygons by J.-P.~Ramis, J.~Sauloy and C.~Zhang \cite{RSZ}. 
Recently, C.~Zhang gave connection formulae for some confluent type basic hypergeometric series \cite{Z1,Z2}.  Zhang also gives the connection formula for the divergent series ${}_2\varphi_0(a,b;-;q,x)$ in \cite{Z0,Z2} where he uses the $q$-Borel-Laplace transformations. In \cite{M0,M1}, the author also gave the connection formula for the Hahn-Exton $q$-Bessel function and the $q$-confluent type function by the using of another kind of the $q$-Borel-Laplace transformations. These resummation methods are powerful tools for connection problems with irregular singular points. 

In this paper, we apply the $q$-Borel-Laplace transformations \textit{for the bilateral series} to the divergent bilateral basic hypergeometric series \eqref{div}. 

\bigskip
\noindent
\begin{defi}
We assume that $f(x)$ is a formal power series $f(x)=\sum_{n\in\mathbb{Z}}a_nx^n$, $a_0=1$.
\begin{enumerate}
\item The $q$-Borel transformation is
\[\left(\mathcal{B}_q^+f\right)(\xi ):=\sum_{n\in\mathbb{Z}}a_nq^{\frac{n(n-1)}{2}}\xi^n\left(=:\psi (\xi )\right).\]
\item For any analytic function $\psi (\xi )$ around $\xi =0$, the $q$-Laplace transformation is
\[\left(\mathcal{L}_{q, \lambda}^+\psi\right)(x):=
\frac{1}{1-q}\int_0^{\lambda\infty}\frac{\varphi (\xi )}{\theta_q\left(\frac{\xi}{x}\right)}\frac{d_q\xi}{\xi}=\sum_{n\in\mathbb{Z}}\frac{\varphi (\lambda q^n)}{\theta_q\left(\frac{\lambda q^n}{x}\right)}.\]
Here, this transformation is given by Jackson's $q$-integral \cite[page 23]{GR}. 
\end{enumerate}
\end{defi}
The definition is a special case of one of the $q$-Laplace transformations in \cite{ZandD, Z0}. The $q$-Borel transformation is the formal inverse of the $q$-Laplace transformation as follows:
\begin{lem}[Zhang, \cite{Z0}]
For any entire function $f(x)$, we have
\[\mathcal{L}_{q,\lambda}^+\circ\mathcal{B}_q^+f=f.\]
\end{lem}
The applications of these transformations can be found in \cite{Z2,Mbi} . We remark that these examples are connection formulae for the bilateral solution of \textit{the first order} $q$-difference equations. But other formulae, especially more higher order and the degenerated (i.e., a confluent) case have not known. In the last section, we give the proof of the main theorem by the using of the $q$-Borel-Laplace transformations and Slater's formula \cite{Slater}.

\section{Basic notations}
In this section, we fix our notations. The $q$-shifted operator $\sigma_q$ is given by $\sigma_qf(x)=f(qx)$. For any fixed $\lambda\in\mathbb{C}^*\setminus q^{\mathbb{Z}}$, the set $[\lambda ;q]$-spiral is $[\lambda ;q]:=\lambda q^{\mathbb{Z}}=\{\lambda q^k;k\in\mathbb{Z}\}$. 
The function $(a;q)_n$ is the $q$-shifted factorial;
\[(a;q)_n:=
\begin{cases}
1, &n=0, \\
(1-a)(1-aq)\dots (1-aq^{n-1}), &n\ge 1,\\
[(1-aq^{-1})(1-aq^{-2})\dots (1-aq^n)]^{-1}, &n\le -1
\end{cases}
\]
moreover, $(a;q)_\infty :=\lim_{n\to \infty}(a;q)_n$ and 
\[(a_1,a_2,\dots ,a_m;q)_\infty:=(a_1;q)_\infty (a_2;q)_\infty \dots (a_m;q)_\infty.\]

\noindent
The basic hypergeometric series with the base $q$ \cite[page 4]{GR} is
\begin{align*}
{}_r\varphi_s(a_1,\dots ,a_r&;b_1,\dots ,b_s;q,x)\\
&:=\sum_{n\ge 0}\frac{(a_1,\dots ,a_r;q)_n}{(b_1,\dots ,b_s;q)_n(q;q)_n}\left\{(-1)^nq^{\frac{n(n-1)}{2}}\right\}^{1+s-r}x^n.
\end{align*}
The radius of convergence is $\infty , 1$ or $0$ according to whether $r-s<1, r-s=1$ or $r-s>1$. 

\noindent
The bilateral basic hypergeometric series with the base $q$ \cite[page 137]{GR} is
\begin{align*}
{}_r\psi_s(a_1,\dots ,a_r&;b_1,\dots ,b_s;q,x)\\
&:=\sum_{n\in\mathbb{Z}}\frac{(a_1,\dots ,a_r;q)_n}{(b_1,\dots ,b_s;q)_n}\left\{(-1)^nq^{\frac{n(n-1)}{2}}\right\}^{s-r}x^n.
\end{align*}

The series ${}_r\psi_s(a_1,\dots ,a_r;b_1,\dots ,b_s;q,x)$ converges on:

\[
  \begin{tabular}{ll}
    $r<s$& $|x|>R:=\displaystyle\left|\frac{b_1b_2\cdots b_s}{a_1a_2\cdots a_r}\right|$ \\
    $r=s$& $R<|x|<1$ \\
     $s<r$ & \textrm{divergent around the origin.}
   \end{tabular}
   \]
\noindent
The theta function of Jacobi is important in connection problems. The theta function of Jacobi with the base $q$ is
\[\theta_q(x):=\sum_{n\in\mathbb{Z}}q^{\frac{n(n-1)}{2}}x^n,\qquad \forall x\in\mathbb{C}^*.\]
The theta function has the triple product identity
\begin{equation}
\theta_q(x)=\left(q,-x,-\frac{q}{x};q\right)_\infty .
\label{triple}
\end{equation}
The theta function satisfies the first order $q$-difference equation
$\theta_q(q^kx)=q^{-\frac{n(n-1)}{2}}x^{-k}\theta_q(x)$, $\forall k\in\mathbb{Z}$. The theta function also has the inversion formula $\theta_q\left(1/x\right)=\theta_q(x)/x$.

We remark that the function $\theta (-\lambda x)/\theta (\lambda x)$, $\forall\lambda\in\mathbb{C}^*$ satisfies a $q$-difference equation
\[u(qx)=-u(x),\]
which is also satisfied by the function $u(x)=e^{\pi i\left(\frac{\log x}{\log q}\right)}$. 

\section{Main theorem}
In this section, we give the proof of the main theorem:
\begin{thm}
For any $x\in\mathbb{C}^*\setminus [-\lambda ;q]$, we have 
\begin{align*}
&\left(\mathcal{L}_{q,\lambda}^+\circ\mathcal{B}^+_q{}_2\psi_1(a_1,a_2;b_1;q,x)\right)(x)\\
=&\frac{(1/a_2,qa_1/a_2,b_1/a_1,q;q)_\infty}{(b_1,q/a_1,a_1/a_2,qa_2/a_1;q)_\infty}\frac{\theta (a_1\lambda /q)}{\theta (\lambda /q)}
\frac{\theta (a_1qx/\lambda )}{\theta (qx/\lambda )}\frac{\theta (x)}{\theta (a_1x)}v_1(x)\notag\\
+&\frac{(1/a_1,qa_2/a_1,b_1/a_2,q;q)_\infty}{(b_1,q/a_2,a_2/a_1,qa_1/a_2;q)_\infty}\frac{\theta (a_2\lambda /q)}{\theta (\lambda /q)}
\frac{\theta (a_2qx/\lambda )}{\theta (qx/\lambda )}\frac{\theta (x)}{\theta (a_2x)}v_2(x).\notag
\end{align*}\end{thm}
Here, $v_1(x)$ and $v_2(x)$ are a fundamental system of unilateral solutions (of equation \eqref{diveq})around  infinity:

\begin{align*}v_1(x)&=\frac{\theta (a_1x)}{\theta (x)}\sum_{n\ge 0}\frac{(qa_1/b_1;q)_n (b_1/a_1a_2x)^n}{(qa_1/a_2;q)_n(q;q)_n},\\
v_2(x)&=\frac{\theta (a_2x)}{\theta (x)}\sum_{n\ge 0}\frac{(qa_2/b_1;q)_n(b_1/a_1a_2x)^n}{(qa_2/a_1;q)_n(q;q)_n}.
\end{align*} 
In the proof of main theorem, Slater's formula for the bilateral series\cite{Slater} plays an important role. In subsection \ref{sl},
we review Slater's formula for a bilateral basic hypergeometric series.
\subsection{Slater's theorem}\label{sl}
Slater gave the following connection formula between the bilateral series ${}_r\psi_{r}(a_1,\dots ,a_r;b_1,\dots b_r;q,x)$ around the origin and the basic hypergeometric function ${}_r\varphi_{r-1}$

\begin{thm}[Slater, \cite{Slater}]\label{sla}
For any $|b_1\cdots b_r/a_1\dots a_r|<|x|<1$, we have
\begin{align*}
&\frac{(b_1,\dots ,b_r,q/a_1,\dots ,q/a_r, x,q/x;q)_\infty}{(qa_1,\dots qa_r,1/a_1,\dots ,1/a_r;q)_\infty}
{}_r\psi_r(a_1,\dots ,a_r; b_1,\dots ,b_r;q,x)\\
&=\frac{a_1^{r-1}(q,qa_1/a_2,\dots qa_1/a_r,b_1/a_1,\dots ,b_r/a_1,a_1x,q/a_1x;q)_\infty}
{(qa_1,1/a_1,a_1/a_2,\dots , a_1/a_r,qa_2/a_1,\dots ,qa_r/a_1;q)_\infty}\\
&\times {}_r\varphi_{r-1}\left(qa_1/b_1,\dots , qa_1/b_r; qa_1/a_2,\dots ,qa_1/a_r;q,\frac{b_1\cdots b_r}{a_1\cdots a_rx}\right)\\
&+\idem (a_1;a_2,\dots ,a_r).
\end{align*}
The notation $\idem (a_1;a_2,\dots ,a_r)$ after an expression stands for the sum of the $r$ expressions obtained from the preceding expression by interchanging $a_1$ with each $a_k$, $k=2,3,\dots , r$.
\end{thm}
A special case of Slater's formula gives Ramanujan's summation formula.
\begin{rem}If we put $r=1$ in theorem \ref{sla}, we obtain Ramanujan's sum for ${}_1\psi_1(a;b;q,x)$\cite[page 57]{Ramanujan}
\begin{align*}
 {}_1\psi_1(a;b;q,z)&=\frac{(q,b/a,az,q/az;q)_\infty}{(b,q/a,z,b/az;q)_\infty}\\
 &=\frac{(b/a,q;q)_\infty}{(b,q/a;q)_\infty}\frac{\theta (-az)}{\theta (-z)}{}_1\varphi_0\left(a;-;q,\frac{q}{az}\right).
\notag
\end{align*}
\end{rem}

\bigskip
We put $r=2$ and take the limit $b_2\to 0$ in theorem \ref{sla},  we obtain the following corollary:

\begin{cor}\label{k1}For any $0<|x|<1$, we have
\begin{align*}
{}_2\psi_2(a_1,a_2;b_1,0;q,x)&=
\frac{(qa_1,qa_2,1/a_2,qa_1/a_2,b_1/a_1,q;q)_\infty}{(b_1,q/a_1,q/a_2,qa_1, a_1/a_2, qa_2/a_1;q)_\infty}\\
&\times \frac{\theta\left(-\frac{a_1x}{q}\right)}{\theta\left(-\frac{x}{q}\right)} 
 {}_1\varphi_1\left(\frac{qa_1}{b_1};\frac{qa_1}{a_2};q, \frac{qb_1}{a_2x} \right)+\idem (a_1;a_2).
\end{align*}
\end{cor}

\subsection{Proof of main theorem}
In this subsection, we give the proof of the main theorem by the using of the $q$-Borel-Laplace transformations.

\begin{proof}We apply the $q$-Borel transformation to the divergent series \eqref{div}. Then, we obtain the following expression for the Borel transform of \eqref{div} by the using of corollary \ref{k1}:
\begin{align*}&\left(\mathcal{B}_q^+{}_2\psi_1(a_1,a_2;b_1;q,x)\right)(\xi )={}_2\psi_2(a_1,a_2;b_1,0;q,-\xi )\\
&=\frac{(qa_1,qa_2,1/a_2,qa_1/a_2,b_1/a_1,q;q)_\infty}{(b_1,q/a_1,q/a_2,qa_1, a_1/a_2, qa_2/a_1;q)_\infty}
\frac{\theta\left(\frac{a_1\xi }{q}\right)}{\theta\left(\frac{\xi}{q}\right)} \\
&\times {}_1\varphi_1\left(\frac{qa_1}{b_1};\frac{qa_1}{a_2};q, -\frac{qb_1}{a_2\xi} \right)+\idem (a_1;a_2)\\
&=:\psi (\xi ).
\end{align*}
We also apply the $q$-Laplace transformation to the function $\psi (\xi )$ as follows:
\begin{align*}
&\left(\mathcal{L}_{q,\lambda}^+\psi(\xi )\right)(x)=\sum_{n\in\mathbb{Z}}\frac{\psi (\lambda q^n)}{\theta (\lambda q^n/x)}\\
&=\frac{(qa_1,qa_2,1/a_2,qa_1/a_2,b_1/a_1,q;q)_\infty}{(b_1,q/a_1,q/a_2,qa_1, a_1/a_2, qa_2/a_1;q)_\infty}\\
&\times\sum_{n\in\mathbb{Z}}\frac{\left(\frac{\lambda}{x}\right)^nq^{\frac{n(n-1)}{2}}}{\theta (\frac{\lambda}{x})} \frac{\theta\left(\frac{a_1\lambda}{q} q^n\right)}{\theta\left(\frac{\lambda}{q}q^n\right)}
\sum_{k\ge 0}\frac{(qa_1/b_1;q)_k \left\{(-1)^kq^{\frac{k(k-1)}{2}}\right\}}{(qa_1/a_2)_k(q;q)_k}\left(-\frac{qb_1}{a_2\lambda}q^{-n}\right)^k\\
&+\idem (a_1;a_2)\\
&=\frac{(qa_1,qa_2,1/a_2,qa_1/a_2,b_1/a_1,q;q)_\infty}{(b_1,q/a_1,q/a_2,qa_1, a_1/a_2, qa_2/a_1;q)_\infty}
\frac{\theta (a_1\lambda /q)}{\theta (\lambda /q)}\frac{\theta (\lambda /a_1x)}{\theta (\lambda /qx)}\\
&\times {}_2\varphi_1(qa_1/b_1,0;qa_1/a_2;q,b_1/a_1a_2x)+\idem (a_1;a_2).
\end{align*}
Therefore, we obtain the conclusion.
\end{proof}
\begin{rem}We set the functions
\[C_1(x):=\frac{(1/a_2,qa_1/a_2,b_1/a_1,q;q)_\infty}{(b_1,q/a_1,a_1/a_2,qa_2/a_1;q)_\infty}\frac{\theta (a_1\lambda /q)}{\theta (\lambda /q)}
\frac{\theta (a_1qx/\lambda )}{\theta (qx/\lambda )}\frac{\theta (x)}{\theta (a_1x)}\]
and 
\[C_2(x):=\frac{(1/a_1,qa_2/a_1,b_1/a_2,q;q)_\infty}{(b_1,q/a_2,a_2/a_1,qa_1/a_2;q)_\infty}\frac{\theta (a_2\lambda /q)}{\theta (\lambda /q)}
\frac{\theta (a_2qx/\lambda )}{\theta (qx/\lambda )}\frac{\theta (x)}{\theta (a_2x)},\]
the new connection formula can be rewritten in the following form:
\[\left(\mathcal{L}_{q,\lambda}^+\circ\mathcal{B}^+_q{}_2\psi_1(a_1,a_2;b_1;q,x)\right)(x)
=C_1(x)v_1(x)+C_2(x)v_2(x).\]
These connection coefficients $C_1(x)$ and $C_2(x)$ are $q$-elliptic functions.
\end{rem}

\section*{Acknowledgements}
The author would like to give heartful thanks to Professor Yousuke Ohyama who provided carefully considered feedback and many valuable comments.

\end{document}